\documentclass[12pt]{article}
\usepackage{amscd,amsfonts,amssymb,multicol,mathtext,color}
\usepackage[tbtags]{amsmath}
\usepackage[cp1251]{inputenc}
\usepackage[english,russian]{babel}
\newtheorem{theorem}{\hskip 0.5 cm Theorem}[section]
\newtheorem{lemma}{\hskip 0.5 cm Lemma}[section]
\newtheorem{note}{\hskip 0.5 cm Remark}[section]
\newtheorem{statement}{\hskip 0.5 cm Proposition}[section]
\newenvironment{proof}{\par\noindent\textbf{Proof.}}{\hfill$\square$}
\newtheorem{definition}{Definition}[section]

\title{Constructive Solutions to the Riemann--Hilbert Problem and Middle Convolution}
\author{Yulia Bibilo\footnote{Corresponding author}\;\footnote{
Institute for Information Transmission Problems, Russian Academy of Sciences,
Bolshoy Karetny per. 19, Moscow, 127994, Russia. Email: y.bibilo@gmail.com} \
  and Galina Filipuk\footnote{Department  of Mathematics, Informatics and Mechanics,
University of Warsaw, Banacha 2, 02-097, Warsaw,
Poland. Email: filipuk@mimuw.edu.pl}}
\date{}

\begin{document}

\maketitle

{\bf Abstract:}
In this paper we present a general scheme to generate constructive solutions to the Riemann--Hilbert problem via middle convolution and illustrate this approach for a  Fuchsian system with four singular points.



\section{Introduction}

A Fuchsian system with a given  number of poles and a given  monodromy representation does not always exist.  There have been a number of studies, discussed below in more details, to obtain necessary and sufficient conditions for such system to exist. In this paper we present a new approach to get explicit solutions of this problem via middle convolution. 
  Middle convolution is related to the Euler transformation of solutions of the Fuchsian systems.

The paper is organized as follows. We shall first  review definitions of the Fuchsian system and its monodromy re\-pre\-sen\-ta\-tion, formulate the Riemann--Hilbert problem and present some of the known methods of  its constructive solutions (other methods and further references can be found in, for example, a survey paper \cite{Gon}).  Then we shall briefly describe an algorithm of middle convolution following \cite{DR, DR2} and, finally, we shall present our new approach (the general scheme) to constructive solutions of the Riemann--Hilbert problem via middle convolution with two illustrative examples.

\subsection{Fuchsian systems}

Let us consider a $(p \times p)$ linear  system of differential equations on the Riemann sphere $\bar{\mathbb{C}}$
\begin{equation}\label{syst_merom}
\frac{dy}{dz}=A(z)y,
\end{equation}
where $A(z)$ is a rational function.
\begin{definition}
Singular points of system $(\ref{syst_merom})$ are the poles of the coefficient matrix $A(z)$. A singular point $a$ is called regular if the fundamental matrix solution of   system $(\ref{syst_merom})$ has at most power growth in every sector  (with the vertex in $a$) of an opening  less then $2\pi$. Otherwise, the singularity  $a$ is called irregular.
\end{definition}
\begin{definition}
If the coefficient matrix $A(z)$ of  system $(\ref{syst_merom})$ has a first order pole at $a$, then the singularity $a$ is called Fuchsian.
\end{definition}
It is known $\cite{Bol3}$ that the Fuchsian singularities are regular. A system with all Fuchsian singularities is called a Fuchsian differential system.
For a  differential system
\begin{equation}\label{syst}
\frac{dy}{dz}=\left(\sum_{i=1}^n \frac{A_i}{z-a_i}\right)y,\quad \sum_{i=1}^n A_i=-A_{n+1} 
\end{equation}
the points $a_1,\ldots,a_n,\,a_{n+1}=\infty$ are the Fuchsian singularities.

\begin{definition}
A Fuchsian singularity $a_i$ of  system $(\ref{syst})$ is called non-resonant if the corres\-ponding matrix $A_i$ has no eigenvalues which differ by a natural number. Otherwise the point  $a_i$ is called a resonant singularity.
\end{definition}

For   system $(\ref{syst})$ one can define its monodromy representation. Consider a fundamental matrix solution  $Y(z)$ of  system $(\ref{syst})$ and its analytic continuation $\tilde{Y}(z)$
along a loop $\gamma$. Then the fundamental solutions are related by $Y(z)=\tilde{Y}(z)G_{\gamma}$,  where the matrix $G_{\gamma}$  is non-degenerate, does not dependent on  $z$ and  depends only on the homotopy class
$[\gamma]$ of the loop $\gamma$. The mapping  $\gamma \rightarrow G_{\gamma}$  defines a representation of the fundamental group of the Riemann sphere with removed singularities  into the group of non-degenerate matrices, i.e., 
\begin{equation}\label{monodr}
\chi:\pi_1 (\bar{\mathbb{C}}\backslash \{a_1,\ldots,a_n,a_{n+1}\},z_0) \rightarrow \mathbb{GL}(p,\mathbb{C}).
\end{equation}
The  generators of the fundamental group  $\gamma_1,\ldots,\gamma_n,\gamma_{n+1}$ (simple loops encircling the  corres\-pon\-ding singular points $\{a_1,\ldots,a_n,a_{n+1}\}$) satisfy   $\gamma_1 \cdot\ldots \cdot \gamma_n=\gamma_{n+1}^{-1}$. 
Let $G_i$ be a matrix equal to $\chi([\gamma_i])$. The matrices $G_1,\ldots,G_n,G_{n+1}$ are called monodromy matrices. The set of monodromy matrices generate a  matrix group with a single group relation $G_1\cdot \ldots \cdot G_n=G_{n+1}^{-1}$.

If one takes another  fundamental matrix solution $X(z)$ of  system $(\ref{syst})$ different from $Y(z)$, then $X(z)=Y(z)C$, where $C$ is a  constant non-degenerate matrix. After  an analytic continuation of $X(z)$   along a loop $\gamma$ one gets a new fundamental solution $\tilde{X}(z)$ given by $\tilde{X}(z)=X(z)\tilde{G}_{\gamma}$. Then the monodromy matrices $G_{\gamma}$ and $\tilde{G}_{\gamma}$ satisfy 
\begin{equation*}
\tilde{G}_{\gamma}=C^{-1} G_{\gamma} C.
\end{equation*}
Hence, the Fuchsian system $(\ref{syst})$ does not only define a unique monodromy representation $\chi$, but it defines a conjugacy class (by a constant matrix) of the monodromy representations.

The monodromy matrices corresponding to non-resonant singularities can be obtained much simpler than that of the resonant singularities.
\begin{lemma}\label{lem_nonres}$\cite[Lect.~6]{Bol3}$
If the point $a_i$ is a non-resonant Fuchsian singularity of  system $(\ref{syst})$, then $G_i\sim e^{2 \pi \mathbf{i} A_i}$.
\end{lemma}

\subsection{The Riemann--Hilbert problem}

The so-called Riemann--Hilbert problem was first mentioned by Riemann at the end of 1850s. In 1900 Hilbert included it in his famous list of mathematical problems under number 21. It is formulated as follows: show that there always exists a linear differential equation of a Fuchsian type (or a Fuchsian system) with given poles and a given monodromy group (see $\cite{An_Bol,Bol3,Bol1,Er88}$ for the details).  

If one compares dimensions of  spaces of parameters (of monodromy matrices and  coefficients of a scalar equation), then it becomes clear that the problem should be stated only for Fuchsian systems and it has a negative solution for scalar Fuchsian equations in general. 
However, in the case of Fuchsian systems it does not always have a positive solution \cite{An_Bol,Bol3,Bol1}. Bolibruch constructed counterexample to the Riemann--Hilbert in 1989 \cite{An_Bol}. Therefore, results containing  necessary or sufficient conditions of the existence of a Fuchsian system with given poles and a given  monodromy representation are of particular importance nowadays. 

Most of the known results are obtained by non-constructive methods. A review of such results can be found in $\cite{Gon}$. However, there also exist constructive approaches. There are several methods to construct a solution to the Riemann--Hilbert problem (i.~e., to construct a Fuchsian differential system  with given singular points and a given monodromy group). They include  
the methods using matrix series and investigating their convergence and analytic continuation (Lappo--Danilevsky's method $\cite{Lap,Lap1,Lap2,Lap3}$, Krylov's method $\cite{Kr}$, Erugin's method $\cite{Er,Er1,Er2,Er3}$ and others) and 
the methods coming from  the theory of isomonodromic deformations (see, for instance, the papers by Boalch $\cite{B}$, Korotkin \cite{Kor} and others).
We shall overview some of the methods in the next section.

\subsection{Constructive solutions}

\subsubsection{Lappo--Danilevsky's method}
The Lappo--Danilevsky (Lappo--Danilevskij) method \cite{Lap,Lap1,Lap2,Lap3} allows one to express the matrices $A_1,\ldots,$ $A_n$ of  system $(\ref{syst})$ in the form of  certain  matrix series depending on either  exponents of monodromy $W_1,\ldots,W_n$ defined below or the monodromy matrices $G_1,\ldots,G_n$ with the requirement that the matrices $W_1,\ldots,$ $W_n$  are close  to zero matrices or the monodromy matrices $G_1,\ldots,G_n$ are close to the identity matrices (with respect to a matrix norm). 

The exponents of monodromy of  system $(\ref{syst})$ are defined as follows. In the neighborhood of each singularity  $a_i$ the  fundamental matrix solution of $(\ref{syst})$ can be written as
\begin{equation}
Y(z)=\bar{Y}(z)(z-a_i)^{W_i},
\end{equation}
where  $\bar{Y}(z)$ is a matrix holomorphic at $a_i$, $\bar{Y}(a_i)\neq 0$ and $W_i$ is a constant matrix. One can find the following fact in $\cite[II]{Lap}$ and $\cite[XII, \textsection 1]{Er}$\footnote{This fact can also be proved by using the Levelt basis, see $\cite[Lect.~6]{Bol3}$.}.
\begin{statement}\label{W_A}
The matrices $W_i$ and $A_i$ corresponding to the Fuchsian singularity $a_i$ of  system $(\ref{syst})$ have the same eigenvalues. 
\end{statement}
The eigenvalues of $W_i$ coincide with the so-called exponents of the Fuchsian system. If the singularity $a_i$ is non-resonant, then the matrix $\bar{Y}^{-1}(z)$ is also holomorphic at $a_i$. 


Explicitly, the matrices $A_1,\ldots,A_n$ of system $(\ref{syst})$ can be expressed in the form  
\begin{equation}\label{lappo_series_W}
A_i=\sum_{\nu=1}^{\infty} \sum_{j_1,\ldots,j_{\nu}}^{1,\ldots,n}R_i(a_{j_1},\ldots,a_{j_{\nu}}|b) W_{j_1}\cdot \ldots \cdot W_{j_{\nu}},
\end{equation}
where $R_i(a_{j_1},\ldots,a_{j_m}|b)$ can be determined recursively. Here   the indices $j_1,\ldots,j_{\nu}$ take values $1,\ldots,n$ independently in  the sum  $\sum_{j_1,\ldots,j_{\nu}}^{1,\ldots,n}$  and $b$ is a complex number used to define the fundamental matrix $Y$ by  $Y(b)=I$.  The series  converges when the norms of the matrices $\left\| W_i \right\|$ are close to zero. 

The matrices $W_i$ and the monodrony matrices $G_i$ corresponding to both resonant and non-resonant Fuchsian singularity $a_i$ are related by
\begin{equation}\label{G_W}
W_i=\frac{1}{2\pi \mathbf{i}}\ln G_i, \quad \mathbf{i}=\sqrt{-1}.
\end{equation}
It is possible to construct a solution to  the Riemann--Hilbert problem  in a class of Fuchsian systems with finite non-resonant   singularities (one should choose a branch of $\ln G_i$ in such a way that the real parts of the eigenvalues of $W_i$ are in the interval $[0,1)$). 

Note that Lappo--Danilevsky proved that the following  formula holds for the coefficients of the system:
\begin{equation}\label{lappo_series_G}
A_i=\frac{1}{2\pi \mathbf{i}} \sum_{\nu=1}^{\infty} \sum_{j_1,\ldots,j_{\nu}}^{1,\ldots,n} Q_i(a_{j_1},\ldots ,a_{j_{\nu}}|b) (G_{j_1}-I)\cdot\ldots \cdot(G_{j_{\nu}}-I),
\end{equation}
where the coefficients $Q_i(a_{j_1},\ldots ,a_{j_{\nu}}|b)$ can be found recursively. The series converges  when every $\left\| G_i-I \right\|$ is small.

\subsubsection{The Gauss system}

Lappo--Danilevsky also considered a particular case of  a Fuchsian $(2\times 2)$ system with three singularities 
\begin{equation}\label{gauss_syst}
\frac{dy}{dz}=\left(\frac{A_1}{z-a_1}+\frac{A_2}{z-a_2}\right)y, \quad A_1+A_2=-A_3.
\end{equation}
System $(\ref{gauss_syst})$ is called the Gauss system. Lappo--Danilevsky fully solved the Riemann--Hilbert problem in terms of the series $(\ref{lappo_series_G})$ for arbitrary monodromy matrices (not only for matrices lying in the neighborhood of the identity matrix). Later on  Krylov \cite{Kr} substantially improved this result and described the solutions in terms of the hypergeometric series and studied their  multivaluedness. Krylov gave all possible solutions for the Riemann--Hilbert problem. It follows from his results it is always possible to construct non-resonant Fuchsian system for irreducible monodromy matrices.

\subsubsection{Erugin's method}\label{erugin_method}

Erugin extended  Lappo--Danilevsky's results for $(2\times 2)$ systems. 
He  constructed a solution to  the Riemann--Hilbert problem   in case of four singular points and $(2 \times 2)$ monodromy matrices \cite{Er,Er1,Er2,Er3}.
More precisely, Erugin solved the Riemann--Hilbert problem formulated as follows. Assume that the $(2 \times 2)$ matrices $W_1,W_2,W_3$ and the singularities $a_1,a_2,a_3,a_4=\infty$ are given. Then find  a $(2 \times 2)$ linear differential system
\begin{equation}\label{er_syst}
\frac{dy}{dz}=\left(\frac{A_1}{z-a_1}+\frac{A_2}{z-a_2}+\frac{A_3}{z-a_3}\right)y, \quad A_1+A_2+A_3=-A_4,
\end{equation}
such that the  fundamental matrix solution of the system $Y(z)$ normalized by $Y(b)=I$  is given by $Y(z)=\bar{Y}_i(z)(z-a_i)^{W_i}$ in the neighborhood of each singularity $a_i$. Here $\bar{Y}_i(z)$ is a second order matrix homomorphic in the neighborhood of $a_i$.
Furthermore, Erugin studied properties of $A_i$ as the function of $W_i$.

Due to the formulation of the problem in this way and  formula $(\ref{G_W})$,  it is possible to find a Fuchsian system $(\ref{er_syst})$ with non-resonant finite singularities. Indeed, the eigenvalues of the matrix $W_i$ depend on the choice of a branch of the logarithm in formula $(\ref{G_W})$, so we can choose eigenvalues of $W_i$ such that their real parts belong to $[0,1)$. The eigenvalues of $A_i$ are equal to the eigenvalues of $W_i$ by the Prop. $\ref{W_A}$. Then there is no pair of eigenvalues of $A_i$ with the integer difference. 
\begin{note}\label{note_Erugin_inf}
In $\cite[XII,\textsection 3,(3.11)]{Er}$ it is shown that $W_4=\frac{1}{2\pi \mathbf{i}}G_4$ with the main branch of the logarithm, where $G_4$ is the monodromy matrix at the infinity.
Thus, the resulting Fuchsian system $(\ref{er_syst})$ is non-resonant at infinity as well.
\end{note}


\section{Middle convolution}

Dettweiler and Reiter's algebraic analogue \cite{DR,DR2} of Katz' middle convolution is a certain trans\-for\-mation of Fuchsian systems which preserves an index of  rigidity. Middle convolution is related to the Euler transformation of the Fuchsian system. 
 There have been numerous studies on middle convolution in recent years, including  applications to special functions, extensions to irregular systems and others. In this section we shall explain main definitions and state some of the results from  $\cite{DR,DR2}$. 

The idea of the algebraic construction is as follows. For a given parameter  one defines monodromy or residue matrices of dimension $(np\times np)$ which are partitioned into blocks and have only one  row of blocks consisting of initial monodromy or residue matrices and the parameter. By finding invariant subspaces and reducing the size of the matrices (if the invariant subspaces are non-zero) one gets a new set of monodromy matrices (or a new Fuchsian system with the same singularities but with new residue matrices). Note that the size of the resulting  matrices (or the system) depends on the choice of the parameter.

\subsection{The multiplicative case $MC_{\lambda}$}\label{MC}

Let us consider a tuple of matrices $G=(G_1,\ldots,G_n)$, $G_i\in \mathbb{GL}(p,\mathbb{C})$, and define $G_{n+1}^{-1}=G_1\cdot \ldots \cdot G_n$.

\begin{definition}
A convolution of the tuple $G$ is a tuple $C_{\lambda}(G)=(M_1,\ldots,M_n)$ given by   
\begin{equation}\label{MC_M}
M_k=\left(\begin{array}{ccc ccc c}
I& 0& & \ldots& & 0&\\
 & \ddots& & & & &\\
 & & I & & & &\\
\lambda(G_1-I) &\ldots & \lambda(G_{k-1}-I) & \lambda G_k &G_{k+1}-I& \ldots & G_n-I\\
 && & & I & &\\
& & & & &\ddots &\\
& 0& &\ldots & &0 &I
\end{array}\right)\in \mathbb{C}^{np \times np},
\end{equation}
where $k=1,\ldots,n$ and $I$ is an $(p\times p)$ identity matrix.
\end{definition}

There are invariant subspaces under each $M_1,\ldots,M_n$
\begin{equation}
\mathcal{K}_k=\left(\begin{array}{c}
0\\
\vdots\\
\ker(G_k-I)\\
\vdots\\
0 \end{array}\right) \subset \mathbb{C}^{np},\quad k=1,\ldots,n,
\end{equation}

\begin{equation}
\mathcal{L}=\cap_{k=1}^n \ker(M_k-I)=\ker(M_1\cdot \ldots \cdot M_n -I).
\end{equation}
Let
\begin{equation}
\mathcal{K}=\oplus_{i=1}^n \mathcal{K}_i.
\end{equation}
Note that if $\lambda \neq 1$, then 
\begin{equation}
\mathcal{L}=\left<\left(\begin{array}{c}
G_2\cdot \ldots \cdot G_n v\\
G_3\cdot \ldots \cdot G_n v\\
\cdots\\
v
\end{array}\right)|v \in \ker(\lambda G_1\cdot \ldots \cdot G_n-I)\right>
\end{equation}
and
$\mathcal{K}+ \mathcal{L}=\mathcal{K}\oplus \mathcal{L}$.

\begin{definition}
Fix an isomorphism between $\mathbb{C}^{np}/(\mathcal{K} + \mathcal{L})$ and $\mathbb{C}^m$ for some $m$.
Then a tuple of matrices $MC_{\lambda}(G)=(\tilde{G}_1,\ldots,\tilde{G}_n)$, $\tilde{G}_i \in \mathbb{GL}(m,\mathbb{C})$, is called
a middle convolution of the tuple  $G$, where
$\tilde{G}_k$ is induced by the action of $M_k$ on $\mathbb{C}^m \cong \mathbb{C}^{np}/(\mathcal{K} + \mathcal{L})$.
\end{definition}
Note that $MC_{\lambda}(G)$ depends on the choice of an isomorphism between $\mathbb{C}^m$ and $\mathbb{C}^{np}/(\mathcal{K} + \mathcal{L})$.
We will denote $\dim(MC_{\lambda}(G))=m$.

Let us define the following conditions for any $\tau\in \mathbb{C}\backslash \{0\} $:  
\begin{equation*}
(*) \;\; \;\cap_{j \neq i} \ker(G_j-I)\cap \ker(\tau G_i-I)=0,\quad i=1,\ldots,n, 
\end{equation*}
\begin{equation*}
(**) \;\; \;\dim \left[ \sum_{j \neq i} {\rm Im}(G_j-I)+{\rm Im}(\tau G_i-I) \right]=p,\quad i=1,\ldots,n.
\end{equation*}



\begin{note}\label{note_*}
The conditions $(*),(**)$ are always fulfilled if  $p>1$ and the matrices $G_1,\ldots,G_n$ generate an irreducible subgroup of $\mathbb{GL}(p,\mathbb{C})$; or $p=1$ and at least two elements of $G_1,\ldots,G_n$, $n\geq 2$, are not identity matrices $\cite[\text{Remark 3.1}]{DR3}$.
\end{note}

\begin{lemma}\label{MC_prop}$\cite{DR,DR2}$
Middle convolution $MC_{\lambda}$ satisfies the following properties:\\
1) If $\lambda \neq 1 $, $\lambda \neq 0 $,  then  
\begin{equation}\label{dim_MC}
\dim(MC_{\lambda}(G))=\sum_{k=1}^n {\rm rk}(G_k-I)-p+{\rm rk}(\lambda G_1\cdot \ldots \cdot G_n-I).
\end{equation}
2) If $\lambda_1 \neq 0$, $\lambda_2 \neq 0$, $\lambda=\lambda_1 \cdot \lambda_2$ and conditions $(*)$, $(**)$   hold, then 
\begin{equation}\label{mult_MC}
MC_{\lambda_2}(MC_{\lambda_1}(G))\cong MC_{\lambda}(G).
\end{equation}
Moreover,
$MC_{\lambda}( MC_{\lambda^{-1}}(G)) \cong G$.

3) If $\lambda \neq 0$, conditions $(*)$, $(**)$ hold and  $G$ is irreducible, then  
$MC_{\lambda}(G)$ is irreducible.
\end{lemma}

\subsection{The additive case $mc_{\mu}$}\label{mc}

Suppose we have a tuple $A$ of matrices $A_1,\ldots,A_n,A_{n+1}$, where $A_i\in \mathbb{C}^{p \times p}$ and  $A_{n+1}=-(A_1+\ldots+ A_n)$.

Define
\begin{equation}
B_k=\left(\begin{array}{ccc ccc c}
O& \ldots& O& O& O& \ldots& O\\
\vdots & \ldots& \vdots & \vdots & \vdots & \ldots& \vdots \\
A_1& \ldots& A_{k-1}& A_k+\mu I& A_{k+1}& \ldots& A_n\\
\vdots & \ldots& \vdots & \vdots & \vdots & \ldots& \vdots \\
O& \ldots& O& O& O& \ldots& O
\end{array}\right) \in \mathbb{C}^{np \times np},
\end{equation}
where $O$ denotes a zero $(p\times p)$ matrix. And there are invariant subspaces under each $B_1,\ldots,B_n$
\begin{equation}
l_k=\left(\begin{array}{c}
0\\
\vdots\\
\ker(A_k)\\
\vdots\\
0
\end{array}\right) \subset \mathbb{C}^{np},\quad k=1,\ldots,n,
\end{equation}
\begin{equation}
l=\cap_{k=1}^n \ker(B_k)=\ker(B_1+\ldots+B_n).
\end{equation}
Let
\begin{equation}
\tilde{l}=\oplus_{i=1}^n l_k.
\end{equation}
Note that if $\mu \neq 0$, then 
\begin{equation*}
l=\left<\left(\begin{array}{c}
v\\
\vdots\\
v
\end{array}\right)|v \in \ker(A_1+\ldots+A_n+\mu I)  \right>
\end{equation*}
and $\tilde{l}+l=\tilde{l}\oplus l$.

\begin{definition}
A tuple of matrices $c_{\mu}(A)=(B_1,\ldots,B_n)$ is called a convolution of the tuple  $A$.
\end{definition}
\begin{definition}
Fix an isomorphism between $\mathbb{C}^{np}/(\tilde{l} + l)$ and $\mathbb{C}^m$ for some $m$.
Then a tuple of matrices $$mc_{\mu}(A)=(\tilde{A}_1,\ldots,\tilde{A}_n),\quad \tilde{A}_i \in \mathbb{C}^{m\times m}$$ is called a middle convolution, where
$\tilde{A}_i$ is induced by the action of $B_i$ on $\mathbb{C}^m \cong \mathbb{C}^{np}/(\tilde{l} + l)$.
\end{definition}

\subsection{A connection between $MC_{\lambda}$ and $mc_{\mu}$}\label{MC_mc}

The following statement is given in \cite{DR}.
\begin{lemma}\label{lem_jourdan}
Let $G=(G_1,\ldots,G_n)$, $G_i \in \mathbb{GL}(p,\mathbb{C})$ and conditions $(*)$, $(**)$  hold. Assume that 
$MC_{\lambda}(G)=(\tilde{G}_1,\ldots,\tilde{G}_n)$, $\lambda \neq 1 $ 
, $\lambda \neq 0 $. Then\\
1)  under the action of $MC_{\lambda}$ every Jordan block $J(\alpha,l)$ appearing in the Jordan decomposition of $G_i$ transforms to the  Jordan block $J(\alpha \lambda,l')$ of  the Jordan decomposition of $\tilde{G}_i$, where
\begin{equation*}
l':=\left\{ \begin{array}{cc} 
l,& \alpha \neq 1,\,\lambda^{-1},\\
l-1,& \alpha = 1,\\
l+1,& \alpha=\lambda^{-1}.
\end{array} \right.
\end{equation*}
All other Jordan blocks which appear  in the Jordan decomposition of $\tilde{G}_i$ are the blocks of the form $J(1,1)$.

2)  under the action of $MC_{\lambda}$ every Jordan block $J(\alpha,l)$ appearing in the Jordan decomposition of $G_{n+1}^{-1}$ transforms to the Jordan block $J(\alpha\lambda,l')$ of the Jordan decomposition of $\tilde{G}_{n+1}^{-1}$, where 
\begin{equation*}
l':=\left\{ \begin{array}{cc} 
l,& \alpha \neq 1,\lambda^{-1},\\
l+1,& \alpha = 1,\\
l-1,& \alpha=\lambda^{-1}.
\end{array}\right.
\end{equation*}
All other Jordan blocks which appear  in the Jordan decomposition of $\tilde{G}_{n+1}^{-1}$ are the blocks of the form $J(\lambda,1)$.
\end{lemma}

The following theorem characterizes a  connection between the additive and multiplicative cases of middle convolution. We will use the following notations: $D_A$ is a  Fuchsian system of differential equations with fixed singularities and  with residue matrices given by the tuple $A$, i.e.,  
\begin{equation*}
\frac{dy}{dz}=\left(\sum_{i=1}^n \frac{A_i}{z-a_i}\right)y, \quad A=(A_1,\ldots,A_n),
\end{equation*}
and $Mon(D_A)$ is a  tuple of monodromy matrices of the system $D_A$.

\begin{theorem}\label{main_th_mc}$\cite{DR}$
Let $A=(A_1,\ldots,A_n)$, $A_i\in \mathbb{C}^{p\times p}$,
$Mon(D_A)=(G_1,\ldots,G_n)$, $\mu \in \mathbb{C}\backslash \mathbb{Z}$, $\lambda=e^{2\pi \mathbf{i} \mu}$. Moreover, assume that the following conditions hold:
\\
1) condition (*);\\ 
2) condition (**);\\
3) ${\rm rk}(A_i)={\rm rk}(G_i-I)$;\\
4) ${\rm rk}(A_1+\ldots+A_n+\mu I)={\rm rk}(\lambda G_1\cdot \ldots \cdot G_n -I)$ 
(or, equivalently,  ${\rm rk}(\mu I-A_{n+1})={\rm rk}(\lambda G_{n+1}^{-1}-I)$).\\
Then $Mon(D_{mc_{\mu-1}(A)})=MC_{\lambda}(Mon(D_A))$ (the fundamental matrix solutions defining monodromy matrices should be chosen in a special way).
\end{theorem}
\begin{note}\label{note_main_th}
It follows from Rem. $\ref{note_*}$, that conditions $1),2)$ of the Theorem $\ref{main_th_mc}$ is true when $p>1$ and $G_1,\ldots,G_n$ generate an irreducible subgroup of $\mathbb{GL}(p,\mathbb{C})$.
\end{note}

\begin{note}
Theorem $\ref{main_th_mc}$ also states that $\dim mc_{\mu-1}(A) = \dim MC_{\lambda}(Mon(D_A))$.
\end{note}

The following diagram illustrates Theorem $\ref{main_th_mc}$:
\begin{equation*}
\begin{array}{ccc}
Mon(D_A)=(G_1,\ldots,G_n)&    \xrightarrow{MC_{\lambda}}  & (\tilde{G}_1,\ldots,\tilde{G}_n)=Mon(D_{\tilde{A}})\\
\chi\left\uparrow\rule{0cm}{0.8cm}\right.\quad \left\downarrow\rule{0cm}{0.8cm}\right. RH & & \chi\left\uparrow\rule{0cm}{0.8cm}\right.\quad \left\downarrow\rule{0cm}{0.8cm}\right. RH\\
A=(A_1,\ldots,A_n)& \xrightarrow{mc_{\mu-1}} & (\tilde{A}_1,\ldots,\tilde{A}_n)=\tilde{A}
\end{array}
\end{equation*}

\section{Main results}

\subsection{A general scheme to  generate constructive solutions}\label{general scheme}

If one has some explicit examples of  solutions to the Riemann--Hilbert problem, then  it is possible to extend the set of these constructive solutions by applying $MC$ to  monodromy matrices and $mc$ to  the residue matrices (provided that conditions of  Theorem~$\ref{main_th_mc}$ are fulfilled). In addition, one can also obtain  some sufficient results about the constructive solutions to the Riemann--Hilbert problem.

Suppose that  we are given a tuple of matrices $G$ and we need to construct a Fuchsian differential system $D_A$ (with fixed singularities  $a_1,\ldots,a_n$) such that 
$Mon(D_A)=G$, i.e.,  we need to obtain a constructive solution to the Riemann--Hilbert problem. In some cases it is possible to   reduce the  problem to the known solutions by using the  additive and multiplicative cases of middle convolution.
 We can use the following idea. We can build a new tuple of matrices $\tilde{G}=MC_{\lambda}(G)$ for a given set $G$ and for some proper parameter $\lambda$ such that the Riemann--Hilbert
problem has a constructive solution for the new set $\tilde{G}$. In this way we can find a tuple of matrices $\tilde{A}$ defining the Fuchsian system with given monodromy matrices $\tilde{G}$. 
 Then if conditions of Theorem $\ref{main_th_mc}$  are satisfied for $\tilde{G}$ and $\tilde{A}$, then by choosing   $MC_{\lambda^{-1}}$, the Fuchsian system defined by $mc_{-\mu+1 }(\tilde{A})$ has monodromy $MC_{\lambda^{-1}}(MC_{\lambda}(G))$. 
The tuples of matrices   $MC_{\lambda^{-1}}(MC_{\lambda}(G))$ and $G$ are equal up to a conjugation by a constant non-degenerate matrix by Lemma $\ref{MC_prop}$. This conjugation corresponds to the choice of the fundamental solution defining monodromy matrices. Therefore, $mc_{-\mu+1}(\tilde{A})$ is the required solution to our problem.

We can  illustrate this idea by the following  diagram, where  the question mark means that we want to find a Fuchsian system for a given tuple of monodromy matrices (and singularities), i.e., want to solve the Riemann--Hilbert problem constructively:
\begin{equation*}
\begin{array}{ccc}
Mon(D_A)=(G_1,\ldots,G_n) & \xrightarrow{MC_{\lambda}}  &( \tilde{G}_1,\ldots,\tilde{G}_n)=Mon(D_{\tilde{A}})\\
& \xleftarrow{MC_{\lambda^{-1}}} &\\
\chi\left\uparrow\rule{0cm}{0.8cm}\right.\quad \left\downarrow\rule{0cm}{0.8cm}\right. RH& & \chi\left\uparrow\rule{0cm}{0.8cm}\right.\quad \left\downarrow\rule{0cm}{0.8cm}\right. RH\\ 
 & \xrightarrow{mc_{\mu-1}}  &
\\? & \xleftarrow{mc_{-\mu+1}} & (\tilde{A}_1,\ldots,\tilde{A}_n)=\tilde{A} .
\end{array}
\end{equation*}

In other words,
\begin{equation*}
\begin{split}
&Mon(D_A)=Mon(D_{mc_{-\mu+1}(\tilde{A})})=MC_{\lambda^{-1}}(Mon(D_{\tilde{A}}))=MC_{\lambda^{-1}}(\tilde{G})=\\
&=MC_{\lambda^{-1}}(MC_{\lambda}(G))=(CG_1C^{-1},\ldots,CG_n C^{-1}), \quad C\in \mathbb{GL}(p,\mathbb{C}).
\end{split}
\end{equation*}

\subsection{Extension of a class of monodromy data for which the Riemann--Hilbert problem has a constructive solution via Erugin's method}
 
In this section  we will prove a theorem that allows us  to reduce the solution of the Riemann-Hilbert problem to Erugin's solution.  

\begin{theorem}\label{th_MC_Erugin}
Let $a_1,a_2,a_3,a_4=\infty$  be  four singular points and assume that the matrices $G_1,G_2,G_3,G_4$, $G_k\in \mathbb{GL}(p,\mathbb{C})$, $p>2$, with $G_1\cdot \ldots \cdot G_4=I$, satisfy the following conditions:\\ 
1) the tuple $G=(G_1,G_2,G_3)$ is irreducible;\\
2) there exists   $\lambda \in \mathbb{C}\backslash\{0,1\}$, such  that by formula $(\ref{dim_MC})$
$$\dim MC_{\lambda}(G)=2.$$
 
Then for a given tuple $G$ of monodromy matrices there exists a constructive solution to the Riemann--Hilbert problem.
\end{theorem}
\begin{note}
To construct solutions to the Riemann-Hilbert problem with our scheme, it is important to use method that allows to construct non-resonan Fuchsian system, such as Erugin's method.
\end{note}
\begin{proof}
We note that condition $\dim MC_{\lambda}(G)=2$ means that the matrices should be such that the dimension of the invariant subspaces $\mathcal{K}+\mathcal{L}$ is equal to $4p-2$. 

To prove the theorem, we will use the general scheme of the constructive solution to the Riemann--Hilbert problem discussed above. Let $MC_{\lambda}(G)=\tilde{G}=(\tilde{G}_1,\tilde{G}_2,\tilde{G}_3)$. Then the tuple of matrices $\tilde{G}$ can be realized as the tuple  of monodromy matrices of some Fuchsian system with non-resonant singular points (by Erugin's  method, see section $\ref{erugin_method}$)
\begin{equation}\label{tilde_A_syst}
\frac{dy}{dz}=\left(\frac{\tilde{A}_1}{z-a_1}+\frac{\tilde{A}_2}{z-a_2}+\frac{\tilde{A}_3}{z-a_3}\right)y.
\end{equation}
 It remains to prove that the tuple of monodromy matrices $\tilde{G}$ and the corresponding differential system $(\ref{tilde_A_syst})$ satisfy conditions of Theorem $\ref{main_th_mc}$. Then $Mon(mc_{-\mu+1}(D_{\tilde{A}}))=MC_{\lambda^{-1}}(\tilde{G})$ (for $\mu=\frac{1}{2\pi \mathbf{i}} \log \lambda$) will be true. 

We have $\dim \tilde{G}=\dim MC_{\lambda}(G)=2>1$, so
to prove conditions 1 and 2 of Theorem $\ref{main_th_mc}$ we should prove that the group generated by $\tilde{G}_1,\tilde{G}_2,\tilde{G}_3$  is an irreducible subgroup of $\mathbb{GL}(2,\mathbb{C})$ (see Remark $\ref{note_main_th}$). It is irreducible subgroup due to case 3) in Lemma $\ref{MC_prop}$. 



Let us prove condition 3: ${\rm rk}(\tilde{A}_i)={\rm rk}(\tilde{G}_i-I)$, $i=1,2,3$. System $(\ref{tilde_A_syst})$ is non-resonant in finite singular points, thus Lemma $\ref{lem_nonres}$ can be used. We have chosen the eigenvalues of the matrix $W_i$ (for Erugin's method) such that their real parts belong to the interval $[0,1)$. Then the eigenvalues of the matrix $\tilde{A}_i$ are also such that their real parts belong to the interval $[0,1)$ by Prop. $\ref{W_A}$. Then Jordan forms of $\tilde{A}_i$ and $\tilde{G}_i$ are agreed (they have the same number of Jordan blocks of equal sizes) because of formula $\tilde{G}_i\sim e^{2 \pi \mathbf{i} \tilde{A}_i}$ in Lemma $\ref{lem_nonres}$.
Also the Jordan form of $\tilde{A}_i$ has a Jordan block $J(0,l)$ if and only if the Jordan form of $\tilde{G}_i$ has a block $J(1,l)$. Thus, condition 3 also holds true. 

To check condition 4, i.e., ${\rm rk}(\mu I+\tilde{A}_{4})={\rm rk}(\lambda^{-1} \tilde{G}_{4}^{-1}-I)$, we will use that $\infty$ is a non-resonant singularity of system $(\ref{tilde_A_syst})$ (see Remark $\ref{note_Erugin_inf}$). Then again Jordan forms of $\tilde{A}_4$ and $\tilde{G}_4$ are agreed and the Jordan form of $\tilde{A}_4$ has a Jordan block $J(0,l)$ if and only if the Jordan form of $\tilde{G}_4^{-1}$ has a block $J(1,l)$ because of the formula $\tilde{G}_i\sim e^{2 \pi \mathbf{i} \tilde{A}_i}$ in Lemma $\ref{lem_nonres}$.

 
\end{proof}

\begin{note}
We can use Krylov's method instead of Erugin's method to construct system with three Fuchsian singular points.
So Theorem $\ref{th_MC_Erugin}$ is true for three singular points and three monodromy matrices. 
\end{note}

\subsection{The first illustrative example}

At first we illustrate Theorem $\ref{th_MC_Erugin}$ with the following simple example. Assume that we are given three matrices 
\begin{equation*}
G_1=\left(\begin{array}{ccc}
\mathbf{i} & 1 & 1\\
0 & 1 &0\\
0 & 0 & 1
\end{array}\right), \quad
G_2=\left(\begin{array}{ccc}
1 & 0 & 0\\
- 2 \mathbf{i} & 2\mathbf{i} & 0\\
2 \mathbf{i} & -\mathbf{i} & -2 \mathbf{i}
\end{array}\right),\quad G_3=(G_1 G_2)^{-1},
\end{equation*}
and we need to construct a Fuchsian system with singular points $a_1,a_2,a_3=\infty$ realizing the given monodromy matrices $G_1,G_2,G_3$.

The tuple $(G_1,G_2)$ is irreducible.
We apply $MC_{\mathbf{i}}$ and get
\begin{equation*}
\tilde{G}_1=\left(\begin{array}{cc}
1 & 2\\
0 & -1
\end{array}\right), \quad
\tilde{G}_2=\left(\begin{array}{cc}
-2 & 0\\
1 & 2
\end{array}\right).
\end{equation*}
We could get residue matrices according to Erugin's or Krylov's method but we will use here the well known facts from analytic differential equations theory.

The next lemma comes from linear algebra.
\begin{lemma}$\cite{Kr}$
Every pair of $(2\times 2)$-matrices $\tilde{A}_i,\tilde{A}_h$ can be transformed (conjugating by a constant non-degenerate matrix $C$) to one of the four types 
\begin{equation} \label{I}
 C\tilde{A}_iC^{-1}=\left(\begin{array}{cc}
\chi_i^{(1)} & 0\\
1 & \chi_i^{(2)}
\end{array}\right), \quad
C\tilde{A}_hC^{-1}=\left(\begin{array}{cc}
\chi_h^{(1)} & \nabla \\
0 & \chi_h^{(2)}
\end{array}\right); \tag{$I$}
\end{equation}
\begin{equation}\label{II}
 C\tilde{A}_iC^{-1}=\left(\begin{array}{cc}
\chi_i^{(1)} & 0\\
1 & \chi_i^{(2)}
\end{array}\right), \quad
C\tilde{A}_hC^{-1}=\left(\begin{array}{cc}
\chi_h^{(1)} & 0 \\
0 & \chi_h^{(2)}
\end{array}\right);  \tag{$II$}
\end{equation}
\begin{equation}\label{III}
 C\tilde{A}_iC^{-1}=\left(\begin{array}{cc}
\chi_i^{(1)} & 0\\
0 & \chi_i^{(2)}
\end{array}\right), \quad
C\tilde{A}_hC^{-1}=\left(\begin{array}{cc}
\chi_h^{(1)} & 0 \\
0 & \chi_h^{(2)}
\end{array}\right); \tag{$III$}
\end{equation}
\begin{equation}\label{IV}
 C\tilde{A}_iC^{-1}=\left(\begin{array}{cc}
\chi_i^{(1)} & 0\\
1 & \chi_i^{(2)}
\end{array}\right), \quad
C\tilde{A}_hC^{-1}=\left(\begin{array}{cc}
\chi_h^{(1)} & 0\\
\nabla  & \chi_h^{(2)}
\end{array}\right); \tag{$IV$}
\end{equation}
where $i \neq h, i,h\in\{1,2\}$, $\nabla \in \mathbb{C}$, $\nabla \neq 0$, $\nabla + (\chi_1^{(1)}-\chi_1^{(2)})(\chi_2^{(1)}-\chi_2^{(2)})\neq 0$.
\end{lemma}

Monodromy matrices $\tilde{G}_1,\tilde{G}_2,\tilde{G}_3$ generate an irreducible subgroup. Then the corresponding Fuchsian system has to be irreducible as well. So we are looking for residue matrices of the type $(I)$.

The formula \begin{equation}\label{chi_eta}
\chi_i^{(k)}=\frac{1}{2\pi \mathbf{i}}{\rm Log}(\eta_i^{(k)}), \quad \mathbf{i}=\sqrt{-1},
\end{equation}
follows from Prop. $\ref{W_A}$, where $\eta_i^{(1)},\eta_i^{(2)}$ are the eigenvalues  of the monodromy matrix $\tilde{G}_i$,
$\chi_i^{(1)},\chi_i^{(2)}$ are the eigenvalues of the residue matrix $\tilde{A}_i$.
Brunches of logarithm in $(\ref{chi_eta})$ should be chosen in such a way that eigenvalues satisfy the following relation 
\begin{equation*}
\chi_1^{(1)}+\chi_1^{(2)}+\chi_2^{(1)}+\chi_2^{(2)}+\chi_3^{(1)}+\chi_3^{(2)}=0,
\end{equation*}
since $\tilde{A}_1+\tilde{A}_2+\tilde{A}_3=0$. 

The value $\nabla$ can be expressed in terms of the eigenvalues $\{\chi_i^{(k)}\}$:  
\begin{equation*}
\nabla = (\chi_{1}^{(1)}+\chi_2^{(1)})(\chi_1^{(2)}+\chi_2^{(2)})- \chi_3^{(1)}\chi_3^{(2)}.
\end{equation*}
Indeed
\begin{equation*}
C\tilde{A}_3 C^{-1}=-C(\tilde{A}_1+\tilde{A}_2)C^{-1}=
-\left(\begin{array}{cc}
\chi_1^{(1)}+\chi_2^{(1)} & \nabla\\
1 & \chi_1^{(2)}+\chi_2^{(2)}
\end{array}\right),
\end{equation*}
\begin{equation*}
\det \tilde{A}_3 = \chi_3^{(1)} \chi_3^{(2)}=(\chi_1^{(1)}+\chi_2^{(1)})(\chi_1^{(2)}+\chi_2^{(2)})-\nabla. 
\end{equation*}
So we have (here we take $C=I$)
\begin{equation*}
\tilde{A}_1=\left(\begin{array}{cc}
-\frac{1}{2} & 0\\
1 & 0
\end{array}\right), \quad
\tilde{A}_2=\left(\begin{array}{cc}
\frac{1}{2}+\frac{1}{2\pi \mathbf{i}}{\rm Log}(2) & \frac{1}{9}\\
0 & \frac{1}{2\pi \mathbf{i}}{\rm Log}(2)
\end{array}\right),
\end{equation*}
\begin{equation*}
\tilde{A}_3=-\tilde{A}_1-\tilde{A}_2=\left(\begin{array}{cc}
-\frac{1}{2\pi \mathbf{i}}{\rm Log}(2) & -\frac{1}{9}\\
-1 & -\frac{1}{2\pi \mathbf{i}}{\rm Log}(2)
\end{array}\right).
\end{equation*}
The system 
\begin{equation*}
\frac{dy}{dz}=\frac{\tilde{A}_1}{z-a_1}+\frac{\tilde{A}_2}{z-a_2},\quad a_1 \neq a_2,
\end{equation*}
has only non-resonant singular points and it is irreducible. Thus, it has monodromy $\tilde{G}_1,\tilde{G}_2,\tilde{G}_3$ which follows from the analogous reasoning.

It is evident this solution of the Riemann--Hilbert problem is not unique.

Next we apply $mc_{\frac{3}{4}}$ according to the main idea.
The resulting system for $G_1,G_2,G_3$ is 
$$\frac{dy}{dz}=\frac{A_1}{z-a_1}+\frac{A_2}{z-a_2},$$
\begin{equation*}
A_1=\left(\begin{array}{ccc}
\frac{1}{4} & \frac{1}{2}+\frac{\log(2)}{2 \pi \mathbf{i}} & \frac{1}{9}\\
0 & 0 & 0\\
0 & 0 & 0
\end{array}\right), \quad
A_2=\left(\begin{array}{ccc}
0 & 0 & 0\\
-\frac{1}{2} & \frac{5}{4}+\frac{\log(2)}{2 \pi \mathbf{i}} & \frac{1}{9}\\
1 & 0 & \frac{3}{4}+\frac{\log(2)}{2 \pi \mathbf{i}}
\end{array}\right).
\end{equation*}

\subsection{The second illustrative example}

In the second example of the constructive solution  to the  Riemann--Hilbert problem  we use  middle convolution and Erugin's method $\cite{Er,Er2}$ (see also \cite{Er1, Er3, Er88}).
In this example we use some of the results from  $\cite{Am_Vas}$ in which Erugin's method for four $(2\times 2)$ matrices in a special form is studied. 

Let us fix four singular points $a_1,a_2,a_3,a_4=\infty$. Define  $G=(G_1,G_2,G_3)$, where the  matrices $G_1,\ldots,G_4$, $G_4=(G_1  G_2  G_3)^{-1}$ are  
  given as follows: 
\begin{equation}\label{ex_G}
\begin{split}
&G_1=\left(\begin{array}{ccc}
a& -2\pi \mathbf{i} \tau & 2\pi \mathbf{i} \tau \\ 
0& 1& 0\\
0& 0& 1
\end{array}\right),
\;\;G_2=\left(\begin{array}{ccc}
1& 0& 0\\
 -2\pi \mathbf{i} \kappa & a & 0\\
0& 0& 1
\end{array}\right),\\
&G_3=\left(\begin{array}{ccc}
1& 0& 0\\
0& 1& 0\\
-2\pi \mathbf{i} \kappa & 0 & a
\end{array}\right),\qquad \tau,\kappa \in (0,1),\; a\in \mathbb{C}\backslash\{0,1\}.
\end{split}
\end{equation}
We  cannot apply  the Lappo--Danilevsky method for the matrices $G_1,\ldots,G_4$ since the norms $\left\|G_i-I\right\|$ are not small when $|a|$ is big. 
Notice that  the matrices $G_1,\ldots,G_4$ satisfy conditions of Theorem $\ref{th_MC_Erugin}$.  Indeed, the matrices $G_1,G_2,G_3$ generate  an irreducible subgroup\footnote{One can verify it by direct calculations, i.e., there is no matrix $C$, such that all matrices $CG_i C^{-1}$ have the same $(2\times1)$ or $(1\times 2)$ zero block under the leading diagonal when $a \notin \{0,1\}$, $\kappa \neq 0,\,\tau\neq 0$.} and, as shown  below, 
$$\dim MC_{a^{-1} } G=2.$$  
Hence, by the Erugin method (see $\cite{Er,Er2}$ and Section~$\ref{erugin_method}$) one can find a constructive solution  to the Riemann--Hilbert problem for  these  matrices.

Let us apply the multiplicative version of middle convolution  $MC_{\lambda}$ with the parameter $\lambda=a^{-1}$ to the tuple $G$.
At first we find the matrices $M_1,M_2,M_3$ according to formula $(\ref{MC_M})$. For each matrix $G_1,G_2,G_3$ one has that $\dim \ker(G_i-I)=2$ and,  moreover,  $\dim \ker (M_1M_2M_3-I)=1$.  
Thus, in the result of the application of the algorithm of middle convolution 
 we get a tuple of  $(2\times 2)$ matrices  
\begin{equation}\label{ex_MC_G}
\tilde{G}_1=\left(\begin{array}{cc}
1& 0\\
-2\pi \mathbf{i}\kappa & 1\end{array}\right),
\;\;\tilde{G}_2=\left(\begin{array}{cc}
1 &- 2\pi \mathbf{i}\tau \\
0& 1
\end{array}\right),
\;\;\tilde{G}_3=\left(\begin{array}{cc}
1 & 2\pi \mathbf{i}\tau \\
0& 1
\end{array}\right).
\end{equation}

Next  we use  Erugin's method for  tuple $(\ref{ex_MC_G})$ following \cite{Am_Vas}.
 We first define the matrices $W_i$ using formula $(\ref{G_W})$  choosing a branch of the logarithm to avoid the resonances. We have 
\begin{equation}
W_1=\left(\begin{array}{cc}
0&0\\
-\kappa&0
\end{array}\right),\;\;
W_2=\left(\begin{array}{cc}
0&-\tau\\
0&0
\end{array}\right),\;\;
W_3=\left(\begin{array}{cc}
0&\tau\\
0&0
\end{array}\right) .
\end{equation}
Furthermore,  from   $\tilde{G}_4=(\tilde{G}_1  \tilde{G}_2  \tilde{G}_3)^{-1}$ we get 
\begin{equation}
\tilde{G}_4=\left(\begin{array}{cc}
1& 0\\
2\pi \mathbf{i}\kappa & 1\end{array}\right),\quad
W_4=\left(\begin{array}{cc}
0&0\\
\kappa&0
\end{array}\right).
\end{equation}
Let us  apply  Erugin's construction to the set $W_1,W_2,W_3,W_4$. We remark that in $\cite{Am_Vas}$ the authors study a more general case. In the result we get the system
\begin{equation}\label{ex_mcA}
\frac{dy}{dz}=\left(\sum_{i=1}^3 \frac{\tilde{A}_i}{z-a_i}\right)y, \quad 
\tilde{A}_i=\left(\begin{array}{cc}
-\sqrt{\frac{-\sigma_{i1} \sigma_{i2}}{\tau \kappa}} & \frac{\sigma_{i2}}{-\tau}\\
\frac{\sigma_{i1}}{-\kappa} & \sqrt{\frac{-\sigma_{i1} \sigma_{i2}}{\tau \kappa}}
\end{array}\right),
\end{equation}
where the values $\sigma_{i1},\sigma_{i2}$ are the sums of the power series in $t=\frac{a_3-a_1}{a_3-a_2}$, $|t|<1$ (see $\cite{Am_Vas}$ for details and calculations), 
in particular, 
\begin{equation}\label{sigma_ik}
\sigma_{ik}={\rm tr} (\tilde{A}_i W_k),
\end{equation}
where 
$$
\tilde{W}_k=\tilde{A}_k+{\rm tr}(\tilde{A}_i \tilde{A}_k)+\sum_{\nu=2}^{\infty} \sum_{j_1,\ldots,j_{\nu}}^{1,2,3}P^*_k(a_{j_1},\ldots,a_{j_{\nu}}|b) {\rm tr}(\tilde{A}_{j_1}\cdot \ldots \cdot \tilde{A}_{j_{\nu}})
$$
and  $$P^*_k(a_{j_1},\ldots,a_{j_{\nu}}|b)=\frac{1}{2\pi \mathbf{i}}P_k(a_{j_1},\ldots,a_{j_{\nu}}|b),$$ 
\begin{equation*}
P_k(a_{j_1}|b)=\left\{ \begin{array}{cc}
2\pi \mathbf{i}, & k=j_1,\\
0,& k \neq j_1,
\end{array}\right.
\end{equation*}
\begin{equation*}
P_k(a_{j_1},\ldots,a_{j_\nu}|b)=\frac{(2\pi \mathbf{i})^\nu}{\nu!},\quad j_1=\ldots=j_\nu=k,
\end{equation*}
\begin{equation*}
P_k(a_{j_1},\ldots,a_{j_\nu}|b)=\int_{a_k}^b\left(\frac{P_k(a_{j_1},\ldots,a_{j_{\nu-1}}|b)}{b-a_{j_{\nu}}}-\frac{P_k(a_{j_2},\ldots,a_{j_\nu}|b)}{b-a_{j_1}}\right)db.
\end{equation*}
We also have 
\begin{equation*}
\tilde{A}^2_j=0,\quad \tilde{A}_j \tilde{A}_k=\rho_{jk}I-\tilde{A}_k \tilde{A}_j, \quad \tilde{A}_j \tilde{A}_k \tilde{A}_j=\tilde{A}_j\rho_{jk},
\end{equation*} 
where $\rho_{jk}={\rm tr}(\tilde{A}_j \tilde{A}_k)$. 
Note that  $\rho_{jk}$ can be expressed as certain power series in $t$ \cite{Am_Vas}.

To get a Fuchsian system with the monodromy matrices $G_1,\ldots,G_4$ we use $mc_{-\mu+1}$ with the parameter $\mu=\frac{1}{2\pi \mathbf{i}}\ln a^{-1}$ for  system $(\ref{ex_mcA})$.
In the result we  get a $(3\times 3)$ system
\begin{equation}\label{ex_A}
\frac{dy}{dz}=\left(\sum_{i=1}^3 \frac{A_i}{z-a_i}\right)y,
\end{equation}
\begin{equation}
\begin{split}
&A_1=\left(\begin{array}{ccc}
-\mu+1 & \sqrt{\frac{-\sigma_{11} \sigma_{12}}{\tau \kappa}} \frac{\sigma_{21}}{\sigma_{11}}-\sqrt{\frac{-\sigma_{21} \sigma_{22}}{\tau \kappa}}&
\sqrt{\frac{-\sigma_{11} \sigma_{12}}{\tau \kappa}} \frac{\sigma_{31}}{\sigma_{11}}-\sqrt{\frac{-\sigma_{31} \sigma_{32}}{\tau \kappa}}\\
0&0&0\\
0&0&0
\end{array}\right),\\
&A_2=\left(\begin{array}{ccc}
0&0&0\\
\sqrt{\frac{-\sigma_{21} \sigma_{22}}{\tau \kappa}} \frac{\sigma_{11}}{\sigma_{21}}-\sqrt{\frac{-\sigma_{11} \sigma_{12}}{\tau \kappa}}&-\mu+1&
\sqrt{\frac{-\sigma_{21} \sigma_{22}}{\tau \kappa}} \frac{\sigma_{31}}{\sigma_{21}}-\sqrt{\frac{-\sigma_{31} \sigma_{32}}{\tau \kappa}}\\
0&0&0
\end{array}\right),\\
&A_3=\left(\begin{array}{ccc}
0&0&0\\
0&0&0\\
\sqrt{\frac{-\sigma_{31} \sigma_{32}}{\tau \kappa}} \frac{\sigma_{11}}{\sigma_{31}}-\sqrt{\frac{-\sigma_{11} \sigma_{12}}{\tau \kappa}}&
\sqrt{\frac{-\sigma_{31} \sigma_{32}}{\tau \kappa}} \frac{\sigma_{21}}{\sigma_{31}}-\sqrt{\frac{-\sigma_{21} \sigma_{22}}{\tau \kappa}}&-\mu+1
\end{array}\right).
\end{split}
\end{equation}
System $(\ref{ex_A})$ solves the Riemann--Hilbert problem for monodromy matrices $(\ref{ex_G})$.

It is easy to see that the constructed system $(\ref{ex_A})$ has non-resonant finite singularities.

Note that the general case in \cite{Am_Vas} with matrices 
$$\tilde{A}_i=\left(\begin{array}{cc}
-\eta_i \theta_i & \eta_i^2\\
-\theta_i^2 & \eta_i \theta_i
\end{array}\right),\quad (\eta_1\theta_2-\eta_2\theta_1)^2+(\eta_1\theta_3-\eta_3\theta_1)^2+(\eta_2\theta_3-\eta_3\theta_2)^2=0, \quad \eta_i,\theta_i \in \mathbb{C}$$ 
 can be treated in a  similar way.

\section{Discussion}

In this paper we showed that it is possible to extend conditions of  constructive solutions of the Riemann--Hilbert problem via middle convolution. We illustrated this idea on Erugin's method and the results in \cite{Am_Vas} for the Fuchsian system of order 2 with four singularities. 

Essentially, theorems similar to Th.~\ref{th_MC_Erugin} can be obtained for other methods of constructive solutions to the Riemann--Hilbert problem (see \cite{Gon} for the overview and further references) by using the general scheme proposed in Section~\ref{general scheme} (for any number of singularities of the Fuchsian system). This considerably extends  the class of monodromy data for which for which the Riemann--Hilbert problem has a (constructive) solution.

There are certain analogues of middle convolution for linear systems with irregular singularities. Therefore, the method used in this paper can be extended to such systems and generalized Riemann--Hilbert problem (if the statement similar to Th. $\ref{main_th_mc}$ can be proved for linear systems with irregular singularities and their monodromy data).

\section*{Acknowledgments}

Galina Filipuk acknowledges the support of the Polish NCN Grant 2011/03/B/ST1/00330. 
Yulia Bibilo carried out the research at the IITP RAS with the support of the
grant from the Russian Foundation for Sciences (project No 14-50-00150).
We thank Renat Gontsov and other members of the seminar "Analytic theory of differential equations" at Steklov Mathematical Institute of RAS  for illuminating discussions.


\begin{thebibliography}{99}
\label{bibl}

\bibitem{An_Bol} D.~V.~Anosov and A.~A.~Bolibruch, The Riemann-Hilbert Problem, Aspects
Math., {\bf E22}, Vieweg, Braunschweig, 1994.

\bibitem{Am_Vas}
V.~V.~Amelkin and M.~N.~Vasilevich, {\it Construction of the second-order Fuchsian systems with
nilpotent irreducible residue matrices}, 
Scientific Publications of the State University of Novi Pazar,
Ser. A: Appl. Math. Inform. and Mech. {\bf 5} (1) (2013),   7--15.

\bibitem{B} Ph.~Boalch, {\it Some explicit solutions to the Riemann–Hilbert problem}, in {\it Differential
Equations and Quantum Groups}, IRMA Lect. Math. Theor. Phys.  {\bf 9}, Eur. Math. Soc.,
Zurich (2007),  85--112.

\bibitem{Bol3} A.~A. Bolibruch, {\it Inverse Monodromy Problems in the Analytic Theory of Differential Equations}, MCCME, Moscow, 2009 (in Russian).

\bibitem{Bol1} A. A. Bolibruch, {\it  Differential equations with meromorphic coefficients}, Proc. Steklov Inst. Math. {\bf 272} (2011), 13--43.	

\bibitem{DR3} M.~Dettweiler  and S.~Reiter, {\it An algorithm of Katz and its application to the inverse Galois problem}, J. Symbolic Comput. {\bf 30} (2000), 761--798.

\bibitem{DR} M.~Dettweiler and S.~Reiter, {\it Middle convolution of Fuchsian systems and the construction of rigid differential systems}, Journal of Algebra {\bf 318} (2007), 1--24.

\bibitem{DR2}  M.~Dettweiler and S.~Reiter, {\it Painlev\'e equations and the middle convolution}, Advances in Geometry {\bf 7} (2007),   317--330.

\bibitem{Er} N.~P.~Erugin, {\it The Riemann Problem}, Nauka i Technika, Minsk, 1982 (in Russian).

\bibitem{Er1} N.~P.~Erugin, {\it The Riemann problem. I}, Differencial'nye Uravnenija {\bf 11} (1975), 771--781 (in Russian).

\bibitem{Er2} N.~P.~Erugin, {\it The Riemann problem. II},  Differencial'nye Uravnenija {\bf 12} (1976),   779--799  (in Russian).

\bibitem{Er3} N.~P.~Erugin, {\it The Riemann problem. III. The case $n=2$ and $m=4$},  Differencial'nye Uravnenija {\bf 13} (1977),   238--254  (in Russian).

\bibitem{Er88} N.~P.~Erugin, {\it The Riemann problem}, Differential Equations {\bf 25} (1989), 907--911.

\bibitem{Gon} R.~R.~Gontsov and V.~A.~Poberezhnyi, {\it Various versions of the Riemann-Hilbert problem for linear differential equations}, 
Russ. Math. Surv. {\bf 63} (4)  (2008), 603--639. 


\bibitem{Kor} D.~Korotkin, {\it Solution of matrix Riemann-Hilbert problems with quasi-permutation monodromy matrices}, Math. Ann. {\bf 329} (2004),  335--364. 

\bibitem{Kr} B.~L.~Krylov, {\it Explicit solution of Riemann problem for Gauss system},
 Tr. Kazan Av. Inst.  {\bf 31} (1956),  203--445 (in Russian).

\bibitem{Lap} J.~A.~Lappo--Danilevsky, {\it Application of Matrix Functions to the Theory of Linear Systems of Ordinary Differential Equations}, GITTL, Moscow, 1957 (in Russian).

\bibitem {Lap1} J. A. Lappo--Danilevskij (J. A. Lappo--Danilevsky), {\it M\'emoires sur la th\'eorie des syst\'emes des \'equations diff\'erentielles lin\'eaires. Vol. I}, Travaux Inst. Physico-Math. Stekloff, {\bf  6}, Acad. Sci. USSR, Leningrad, 1934, 1–256.

\bibitem {Lap2}	J. A. Lappo--Danilevskij (J. A. Lappo--Danilevsky), {\it M\'emoires sur la th\'eorie des syst\'emes des \'equations diff\'erentielles lin\'eaires. Vol. II}, Travaux Inst. Physico-Math. Stekloff, {\bf  7}, Acad. Sci. USSR, Moscow–Leningrad, 1935, 5–210.

\bibitem {Lap3} J. A. Lappo--Danilevskij (J. A. Lappo--Danilevsky), {\it M\'emoires sur la th\'eorie des syst\'emes des \'equations diff\'erentielles lin\'eaires. Vol. III}, Travaux Inst. Physico-Math. Stekloff, {\bf  8}, Acad. Sci. USSR, Moscow–Leningrad, 1936, 5–206.

\end{thebibliography}
\end{document}